\documentclass[11pt,fleqn,titlepage]{report}
\usepackage{latexsym,amsmath,amssymb,amsfonts,epsfig}
\usepackage[all]{xy}

\providecommand{\da}{\downarrow}

\providecommand{\CC}{{\mathbb{C}}}
\providecommand{\RR}{{\mathbb{R}}}

\providecommand{\ZZ}{{\mathbb{Z}}}

\providecommand{\OO}{{\mathcal O}}

\providecommand{\QED}{{\hfill{} $\Box$}}

\providecommand{\DD}{{\cal D}} 
\providecommand{\GG}{{\cal G}}
\providecommand{\HH}{{\cal H}}
\providecommand{\KK}{{\cal K}}
\providecommand{\PP}{{\mathbb P}}
\providecommand{\calP}{{\cal P}}
\providecommand{\Sch}{{\cal S}}
\providecommand{\calS}{{\cal S}}
\providecommand{\UU}{{\cal U}}

\providecommand{\Mat}{{\cal M}}  
\providecommand{\THM}{{\mathbb{T}_H M}}
\providecommand{\TM}{{\mathbb{T}M}}
\providecommand{\THG}{{\mathbb{T}_H} G}
\providecommand{\THU}{{\mathbb{T}_H} U}

\providecommand{\PP}{{\mathbb{P}}}

\providecommand{\Lg}{{\mathfrak g}}

\providecommand{\sm}{{\Gamma^\infty}}  
\newcommand{\ang}[1]{\langle #1 \rangle} 

\newtheorem{definition}{Definition}
\newtheorem{lemma}[definition]{Lemma}
\newtheorem{proposition}[definition]{Proposition}
\newtheorem{corollary}[definition]{Corollary}
\newtheorem{theorem}[definition]{Theorem}

\begin{document}

\title{The Atiyah-Singer Index Formula for Subelliptic Operators on Contact Manifolds, Part I}
\author{Erik van Erp \\
The University of Pennsylvania \\
erikerp@math.upenn.edu}

\date{April 2006}
\maketitle

\begin{center}{\Large \bf Abstract}\end{center}
\vskip 12pt

\noindent 
The Atiyah-Singer index theorem gives a topological formula for the index of an elliptic differential operator. The topological index depends on a cohomology class that is constructed from the principal symbol of the operator.
On contact manifolds, the important Fredholm operators are not elliptic, but hypoelliptic. Their symbolic calculus is noncommutative, and is closely related to analysis on the Heisenberg group.

For a hypoelliptic differential operator in the Heisenberg calculus on a contact manifold we construct a symbol class in the K-theory of a noncommutative $C^*$-algebra that is associated to the algebra of symbols. There is a canonical map from this analytic K-theory group to the ordinary cohomology of the manifold, which gives a de Rham class to which the Atiyah-Singer formula can be applied. We prove that the index formula holds for these hypoelliptic operators.

Our methods derive from Connes' tangent groupoid proof of the index theorem.

\tableofcontents

\chapter{Introduction}

In this paper we present the solution to the index problem for a certain class of {\em hypoelliptic} operators on contact manifolds. The operators in question are `elliptic'  (meaning that they have an invertible symbol) in the {\em Heisenberg calculus}, which is a symbolic calculus of (pseudo)differential operators naturally suited to the contact structure. On a compact contact manifold such operators, though not elliptic in the usual sense, are still Fredholm [BG], [CGGP]. 
The problem of deriving a corresponding index formula has been considered extensively by Melrose and Epstein [EM], [Me], [Ep]. 
Their `bottom-up' approach takes as its starting point the index theorem for Toeplitz operators on contact manifolds of Boutet de Monvel [Bo], and derives formulas for progressively larger classes of operators.
Their work clarifies a great deal about the fundamental structures that play a role in the problem. However, their most general formula for the index is quite involved.
This is related to the fact that they only managed to define a $K$-cocycle for a restricted class of hypoelliptic operators (the so-called `Hermite operators', a class that includes the `higher rank' Toeplitz operators.)

Our approach departs, at least initially, radically from that of Melrose and Epstein. We take as our inspiration the abstract {\em tangent groupoid} proof of the classical Atiyah-Singer index theorem suggested by Connes [Co]. Once the relevant definitions have been established in the context of the Heisenberg calculus (a nontrivial matter), this approach results, in one go, in an elegant, but rather abstract, index theorem for the operators under consideration. 
However, while ideas from noncommutative geometry have guided the way to the formulation and proof of this theorem, certain structures appearing in the work of Melrose and Epstein have proved useful in the second stage of our project, when we turn our initial abstract result into a more explicit and computable formula. 

We present our results in two papers, reflecting the two stages by which we arrived at our final theorem.
The current paper is, in essence, a summary of our  Ph.D. thesis [Er].
Here we take the first step by showing that the {\em symbol} of a Heisenberg `elliptic' operator  defines, in a natural way, an element in the $K$-theory of a certain {\em noncommutative $C^*$-algebra}.
If $M$ is a compact contact manifold, with contact hyperplane bundle $H\subset TM$, then the relevant $C^*$-algebra is the convolution algebra of a {\em smooth groupoid} $T_HM$. This groupoid is just the tangent bundle $TM$, but each fiber is equipped with the structure of a Heisenberg group (the so-called {\em osculating groups} introduced by Stein and Folland [FS1]).
For a Heisenberg `elliptic' operator $P$ we construct a $K$-theory element,
\[ [\sigma_H(P)] \in K_0(C^*(T_HM)).\]
We then modify Connes' construction of the tangent groupoid, replacing the tangent bundle $TM$ by the groupoid $T_HM$. This leads to a map in $K$-theory,
\[ K_0(C^*(T_HM)) \to \ZZ. \]
We prove that this map computes the index of $P$. 

Since the fibers of $T_HM$ are nilpotent groups, a well-known result in noncommutative geometry gives a natural isomorphism,
\[ \psi\;\colon\;K_0(C^*(T_HM)) \cong K^0(T^*M) .\] 
This means that we can identify our noncommutative symbol $[\sigma_H(P)]$ 
with an element in the topological $K$-theory group $K^0(T^*M)$.
We end the present paper with a novel groupoid argument to prove that, in fact, the index of $P$ is computed by means of the classical formula of Atiyah-Singer,
\[ {\rm Index}\, P = \int_{T^*M} {\rm Ch}(\psi[\sigma_H(P)])\wedge{\rm Td}(M) .\]
While the isomorphism $\psi$ is natural, it depends on the {\em Connes-Thom Isomorphism} in analytic $K$-theory, and is thus highly nontrivial and not explicitly computable.
In the second stage of our project we show how to directly construct the {\em topological} $K$-theory class associated with $P$, thus bypassing the Connes-Thom Isomorphism in the final formulation of the theorem, and arriving at an explicitly  computable formula.
Here some of the structures elucidated by Melrose and Epstein play a role,
but in the context of analytic $K$-theory they assume new meanings. 

We will see that, in our approach, the Toeplitz index theorem of Boutet de Monvel, rather than serving as a starting point, will be derived as a corollary of the general result.
These more explicit results (not contained in our Ph.D. thesis) are presented in a separate paper [Er2]
  
\vskip 6pt

For the reader who is unfamiliar with the Heisenberg calculus, and for purposes of self-containment, Chapter \ref{chapter:calc} contains a proof of Fredholmness for operators that are `elliptic' in this calculus.
We develop the calculus only for differential operators, which suffices for our purposes here.
The Heisenberg calculus for differential operators is much simpler than the corresponding pseudodifferential theory.
As an introduction to the pseudodifferential theory we recommend [Ep].
A detailed treatment is given in [EM2], which is, to date, unpublished.
The usual references are [BG] (which contains a good overview of the history of the topic) and [Ta1].
An alternative approach (based on operator kernels, rather than symbols) is found in [CGGP].

The reader familiar with the Heisenberg calculus should skip to Chapter \ref{chapter:index},
where we derive our index theorem.
However, since our treatment of the {\em osculating groups} is slightly unconventional, we recommend at least reviewing section \ref{section:calc} for a definition of the groupoid $T_HM$.
The noncommutative $K$-theory class $[\sigma_H(P)]$ is constructed in section \ref{section:Ksymbol}.
Theorem \ref{hypind50:thm:analindx} in section \ref{section:aind} states the index theorem as it emerges from the tangent groupoid machinery, while Theorem \ref{MAIN} in section \ref{section:tind}, the main result of this paper, requires a further and novel groupoid argument.

\chapter{The Heisenberg Calculus on Contact Manifolds}\label{chapter:calc}

\section{Index theory}

We briefly review the basics of index theory for elliptic operators, to highlight some of the features that play a role in subelliptic theory.
Let $P$ be a differential operator acting on smooth sections in vector bundles $E,F$ over a manifold $M$,
\[ P\;\colon\; L^2(E)\to L^2(F) .\]
We treat $P$ as an unbounded closed Hilbert space operator.
Such an operator is {\em Fredholm} if its kernel and cokernel are finite dimensional,
and the {\em Fredholm index} of $P$ is the difference between their dimensions,
\[ {\rm Index}\, P = {\rm dim Ker}\, P - {\rm dim Ker}\, P^* .\] 
If $P$ is elliptic of order $d$, and $M$ closed, then $P$ satisfies {\em elliptic estimates}, which state that for every differential operator $A$ of order $\le d$ there is a constant $C>0$ such that
\[ \|Au\| \le C(\|Pu\| + \|u\| ) ,\]
for every smooth section $u\in \sm(E)$.
By an application of the Rellich Lemma in Sobolev theory, these basic estimates imply that $P$ is Fredholm.

To prove the estimates, choose coordinates on $M$ and local trivializations of $E,F$, and write
\[ P = \sum_{|\alpha|\le d} a_\alpha(x)\partial^\alpha .\]
The behaviour of $P$ near a point $m\in M$ is closely approximated by the {\em homogeneous, constant coefficient} operators
\[ P_m = \sum_{|\alpha| =  d} a_\alpha(m)\partial^\alpha \]
obtained by `freezing the coefficients' of the principal part of $P$ at $m$.
Each operator $P_m$ is invariantly defined as an operator on the tangent space $T_mM$.
Using Fourier theory, one easily derives elliptic estimates for each $P_m$ (where the operator $A$ appearing in the estimates is a constant coefficient operator on $T_mM$ as well). Then, comparing $P$ with $P_m$, one shows that the same estimates hold for $P$ in a small neigborhood of $m$. Compactness of $M$ finishes the proof.

The Fredholm index is a robust invariant, and it only depends on the homotopy type of the principal symbol $\sigma(P)$,
\[ \sigma(m,\xi) = \sum_{|\alpha| =  d} a_\alpha(m)(i\xi)^\alpha ,\]
which is a {\em topological} gadget derived from the family of constant coefficient operators $P_m$.
Its homotopy class defines a cocycle in $K$-theory with compact supports,
\[ [\sigma(P)] \in K^0(T^*M) .\]
The index problem for elliptic operators is to find a topological procedure that derives the Fredholm index of $P$ from its symbol $\sigma(P)$.
Stated more explicitly, we must find a topological characterization of the analytic index map
\[ K^0(T^*M) \to \ZZ \;;\; [\sigma(P)] \mapsto {\rm Index}\, P .\]
In its cohomological form, the Index Theorem of Atiyah and Singer gives the solution as,
\[ {\rm Index}\, P = \int_{T^*M} {\rm Ch}(\sigma(P))\wedge{\rm Td}(M) .\]
It would seem that the Atiyah-Singer formula {\em cannot} apply to operators that are not elliptic, because the Chern character ${\rm Ch}(\sigma(P))$ is only defined when the symbol gives a well-defined element in $K^0(T^*M)$, which is the case precisely when $P$ is elliptic. 
But there is at least one example of non-elliptic Fredholm operators to which the Atiyah-Singer formula applies, namely {\em Toeplitz operators} 
\[ T_a = SaS .\]
Here $S$ denotes the Szeg\"o projector associated with a strictly pseudoconvex domain $X$, and $a$ is a smooth function on the boundary $M=\partial X$. 
The Szeg\"o projector is the orthogonal projection of $L^2(M)$ onto the Hardy space $H^2(M)$, the $L^2$ closure of the space of smooth functions on $M$ that extend holomorphically to $X$.
If $a$ is invertible, then $T_a$ is a Fredholm operator on $H^2(M)$, and Boutet de Monvel ([Bo]) proved that,
\[ {\rm Index}\, P = \int_M {\rm Ch}(a)\wedge{\rm Td}(M) .\]
The invertible function $a$ defines an element in $K^1(M)$, and it plays the same role in the formula as the symbol of an elliptic operator. 
In fact, in the {\em Heisenberg calculus} of pseudodifferential operators on a contact manifold, the function $a$ is, in fact, the symbol of $T_a$, and its invertibility makes $T_a$ an `elliptic' operator in that calculus.
The result presented in this paper accomplishes for {\em differential operators} on $M=\partial X$
what Boutet de Monvel did for Toeplitz operators. For a differential operator $P$ that is `elliptic' in the Heisenberg calculus, we identify the appropriate $K$-theory class
\[ \psi[\sigma_H(P)] \in K^0(T^*M) \]
associated to its Heisenberg symbol, and prove that
\[ {\rm Index}\, P = \int_{T^*M} {\rm Ch}(\psi[\sigma_H(P)])\wedge{\rm Td}(M) .\]
Before turning to the proof of this result, we give a brief and self-contained treatment of the Heisenberg calculus for differential operators.

\section{Contact manifolds}

For the basic facts on contact manifolds see [Ar].

Throughout this paper $M$ denotes a closed (compact, without boundary) contact manifold.
A contact manifold $M$ has odd dimension $2n+1$, and is equiped with a hyperplane bundle $H\subseteq TM$ (i.e., a vector bundle with fiber dimension $2n$),
such that an arbitrary non-vanishing local $1$-form $\theta$ on $M$ with $\theta(H) = 0$ has the property that $\theta(d\theta)^n$ is a nowhere vanishing volume form.
If $\theta$ satisfies the requirements, then so does $f\theta$ for any non-zero function $f$.
In this paper we take the point of view that the contact structure is represented by the contact hyperplane bundle $H$, rather than the contact form $\theta$. 

The canonical example of a contact manifold is the {\em Heisenberg group} $G=\RR^{2n+1}$.
In coordinates $(x,y,t)=(x_1,\ldots,x_n,y_1,\ldots,y_n,t) \in \RR^{2n+1}$ the group operation  is expressed as
\[ (x,y,t)(x',y',t') = (x+x',y+y',t+t'+\frac{1}{2}\sum_{j=1}^n (x_jy'_j - y_jx'_j)) .\]
The right invariant vector fields on $G$ are
\begin{align*}
  X_j & = \frac{\partial}{\partial x_j} + \frac{1}{2} y_j \frac{\partial}{\partial t},\; j=1,\ldots,n, \\
  Y_j & = \frac{\partial}{\partial y_j} - \frac{1}{2} x_j \frac{\partial}{\partial t},\; j=1,\ldots,n, \\
  T & = \frac{\partial}{\partial t}.
\end{align*}
On $G$, let $H\subseteq TG$ be the right invariant hyperplane bundle spanned by the vector fields $X_j$ and  $Y_j$ ($j=1,\ldots,n$). 
Correspondingly, we can take $\theta$ to be the right invariant $1$-form that is equal to $dt$ at the origin.
We have
\[ \theta = dt + \frac{1}{2} \sum_{i=1}^{n} \left( x_j dy_j - y_j dx_j \right) ,\]
and so $d\theta=\sum dx_idy_i$, and $\theta(d\theta)^n$ is the standard volume form on $\RR^{2n+1}$.
The form $\theta$ is the {\em canonical contact form} on the Heisenberg group.

The analysis of subelliptic operators on contact manifolds is greatly simplified by the following fact.

\begin{theorem}
Every contact manifold $(M,H)$ is locally isomorphic to an open subset of the Heisenberg group $G$ with its canonical contact structure.
More precisely, for every $m\in M$ there exists an open set $U\subseteq M$ and a chart $\psi\colon U\to \RR^{2n+1}$ such that the contact form on $M$ is equal to the pull-back of the canonical contact form on $\RR^{2n+1}$. 
\end{theorem}
The proof is an easy application of Darboux's theorem (see [Ar]).

\vskip 6pt
\noindent {\bf Example.}
Contact manifolds arise naturally in the study of boundary value problems in complex analysis.
Let $M=\partial X$ be the smooth boundary of a domain $X\subseteq \CC^{n+1}$.
Following standard notation, let $T^{1,0} \oplus T^{0,1}$ denote the splitting of the complexified tangent bundle of $\CC^{n+1}$ into holomorphic and anti-holomorphic vectors.
Since $M$ has real dimension $2n+1$ it is certainly  not a complex manifold, but it has extra structure which derives from the holomorphic structure on $\CC^{n+1}$, as follows.
Let $T_\CC M$ denote the complexified tangent bundle $TM\otimes \CC$, and
\[ T^{1,0}M = T_\CC M \cap T^{1,0} \;,\; T^{0,1}M= T_\CC M \cap T^{0,1}. \]
Then the complex vector bundle 
\[ H_\CC = T^{1,0}M\oplus T^{0,1}M  \]
is a (complex) codimension 1 subbundle of $T_\CC M$.
Therefore $H = H_\CC \cap TM$  is a {\em real} codimension 1 subbundle of $TM$,
and under suitable conditions this bundle defines a contact structure.

Suppose the region $X$ is given in the form $X=\{z|\rho(z)<0\}$, with $\rho$  a smooth function on $\CC^{n+1}$ whose gradient $d\rho$ is non-zero on the boundary $M = \rho^{-1}(0)$.
Then the $(1,0)$-form $\theta = -i \partial \rho$ restricts to a real $1$-form on $TM$,
and the kernel of $\theta$ (as a form on $TM$) is the hyperplane bundle $H$ defined above.
A hermitian form on $T^{1,0}M$ (the so-called {\em Levi form}) is given by
\[  d\theta =i\partial \bar{\partial}\rho \;\colon\; T^{1,0}M \otimes T^{0,1}M \to \CC.\]
One easily verifies that  $\theta$ is a contact form on $M$ if and only if the Levi form is non-degenerate.
The domain $X$ is called {\em strictly pseudoconvex} if the Levi form is positive definite,
and we see that the boundary of such a domain has a natural contact structure.

It is precisely for such manifolds $M=\partial X$, where $X$ is a strictly pseudoconvex domain, that the index theorem for Toeplitz operators is formulated.

\section{A noncommutative symbolic calculus}\label{section:calc}

The essential feature of subelliptic theory on contact manifolds is that vector fields transversal to the contact hyperplane bundle $H$ are treated as {\em second order} operators.
From this simple idea a calculus is derived quite naturally. 
First, this notion of order of vector fields generates a filtration on the algebra of differential operators $\calP$ on $M$. 
By abstract nonsense, there is a {\em graded} algebra $\calS$ associated to $\calP$, and the {\em principal symbol} of an operator $P\in \calP$ is simply the image of its highest order part in $\calS$.
It turns out that elements in $\calS$ can be identified with sections in a bundle of graded algebras $\UU = \{\UU_m,m\in M\}$.
This means that if we can define the principal part of the operator $P$ {\em at a point} $m\in M$ as an element in a noncommutative algebra $\UU_m$.
Finally, each fiber $\UU_m$ is realized as the universal enveloping algebra $\UU(\Lg_m)$ of a graded nilpotent Lie algebra $\Lg_m$, which has a simple and concrete definition.
This, in essence, is the Heisenberg calculus. We now turn to the details.

We start with a filtration on the algebra of differential operators $\calP$ defined as follows: any vector field $X\in \sm(H)$ defines an order $1$ operator, while any other vector field has order $2$.
Let $\calP^d$ denote the linear span of monomials $\Pi X_k$, where the sum of the orders of the vector fields $X_k$ is less or equal $d$.
We set $\calP^0 = C^\infty(M)$.
A differential operator $P$ has order $d$ in the Heisenberg calculus if $P\in \calP^d\setminus \calP^{d-1}$.
  
Abstractly, the algebra of symbols for the Heisenberg calculus is simply the graded algebra $\calS$ associated to the filtered algebra $\calP$ in the usual manner,
\[  \calS = \bigoplus \calS^d \;;\; \calS^d = \calP^d/\calP^{d-1}.\]
The {\em principal part} of a differential operator $P$ of order $d$ is simply the image of $P$ under the quotient map 
\[ \sigma_H^d\;\colon\; \calP^d \to \calS^d.\]
Clearly, if $P$ has Heisenberg order $k$, and $Q$ has Heisenberg order $l$, then
\[ \sigma_H^k(P)\sigma_H^l(Q) = \sigma_H^{k+l}(PQ). \]
The commutator of a smooth function and a vector field $[X,f]=X.f$ is again a smooth function.
It follows that for any operator $P$ of Heisenberg degree $d$, the commutator $[P,f]$ has Heisenberg degree $d-1$.
Therefore $\calS^0 = C^\infty(M)$ commutes with all of $\calS$.
This crucial fact allows us to localize elements in $\calS$, realizing its elements  as sections in a bundle $\UU$, each of whose fibers $\UU_m$ has the structure of a graded algebra.
We now derive a concrete model for this abstract element
\[ P_m = \sigma_H^d(P)(m) \in \UU_m .\]
Elements in $\UU_m$ are equivalence classes of operators.
Two operators $P,P'$ of order $d$ represent the same class in $\UU_m$ if  
\[ P = P' + \sum_k f_kA_k + \text{lower order terms} ,\]
where $A_k$ are operators of order $d$, and $f_k$ are smooth functions with $f_k(m)=0$.
Suppose we have another element $Q\sim Q'$ in $\UU_m$, with $Q = Q' + \sum_l g_lB_l + \cdots$.
Then
\begin{align*}
 PQ & = P'Q' + P'\sum g_l B_l  + \sum f_k A_k Q' + \sum f_kA_kg_lB_l + \text{lower order terms} \\ 
    & = P'Q' + \sum g_lP'B_l + \sum f_kA_kQ' + \sum f_kg_lA_kB_l + \text{lower order terms} .
\end{align*}
Thus, $PQ$ and $P'Q'$ represent the same element in $\UU_m$,
and we have a well-defined graded algebra structure on $\UU_m$.

Consider the image of the set of vector fields $\sm(TM)\subseteq \calP$ in $\UU_m$.
The principal parts of vector fields are sections in the bundle $H\oplus N\subseteq \calS$, where $N$ denotes the quotient line bundle $TM/H$ (the symbols of vector fields of degree $2$).
Restricting to $m\in M$ we get
\[ H_m\oplus N_m \subseteq \UU_m .\]
Since the commutator of two vector fields is again a vector field,
the subset $H_m\oplus N_m$ is closed under commutators in the algebra $\UU_m$, and has the structure of a {\em graded Lie algebra} ($H_m$ has degree $1$, $N_m$ degree 2).
We call it the {\em osculating Lie algebra}, and denote it by $\Lg_m$.

Since vector fields generate the algebra of differential operators $\calP$,
the subset $\Lg_m$ generates $\UU_m$ as an algebra.
In fact, one easily verifies that $\UU_m$ is naturally isomorphic to the universal enveloping algebra $\UU(\Lg_m)$ of $\Lg_m$ (by the Poincar\'e--Birkhoff--Witt Theorem), and we obtain,
\[ P_m\in \UU(\Lg_m) .\]
Observe that $\Lg_m$ is two-step nilpotent, with $[H_m,H_m]\subseteq N_m$ and $[N_m,\Lg_m]=\{0\}$. 
By definition of the algebra structure on $\UU_m$, the bracket of two elements $Y_1,Y_2\in H_m$ is calculated by choosing two vector fields $X_1,X_2\in \sm(H)$ with $Y_j=X_j(m)$,
and letting
\[ [Y_1,Y_2] = [X_1,X_2](m) \,{\rm mod}\, H_m .\]  
On the Heisenberg group $G$ with its canonical contact structure,
we see that all osculating Lie algebras $\Lg_m$ are isomorphic to the Lie algebra $\Lg$ of $G$.
By Darboux's theorem, the same is true for contact manifolds.

To sum up the discussion, and make it more explicit,
choose vector fields $X_1,\cdots, X_{2n}$ that span the fibers of $H$ at each point in some open set of $M$, and let $X_0$ be any vector field transversal to $H$. 
Such a collection of vector fields is called an {\em $H$-frame}.
A differential operator $P$ of Heisenberg order $d$ can be locally presented by means of these vector fields,
\[ P = \sum_{|\alpha|\le d} a_\alpha X^\alpha .\]
As usual, $\alpha = (\alpha_0,\alpha_1,\ldots,\alpha_{2n})$ is a multi-index, while
the expression $X^\alpha$ is analogous to the usual notation $\partial^\alpha$, and is shorthand for
\[ X^\alpha = X_0^{\alpha_0} X_1^{\alpha_1} \cdots X_{2n}^{\alpha_{2n}} .\]
Observe that the meaning of $X^\alpha$ depends on the order of the vector fields in the $H$-frame.
The Heisenberg degree of the monomial $X^\alpha$ is given by,
\[  |\alpha| = 2\alpha_0 + \alpha_1 + \cdots + \alpha_{2n} .\]
Freezing coefficients of the principal part of $P$ at the point $m\in M$, following the procedure described above, we obtain
\[ P_m = \sum_{|\alpha| =  d} a_\alpha(m) Y^\alpha \in \UU(\Lg_m),\]
with $Y_j = X_j(m) \in \Lg_m = H_m \oplus N_m$.

Let $G_m$ be the simply connected nilpotent Lie group associated with $\Lg_m$.
By the Campbell--Baker--Hausdorff formula, the group operation on $G_m=H_m\oplus N_m$ is given by the {\em quadratic} formula,
\[ u\cdot v  = u+v+\frac{1}{2}[u,v] .\]
This group is called the {\em osculating group} at the point $m$.
We let $T_HM$ denote the bundle $H\oplus N$ with the nilpotent group structure of $G_m$ on its fibers.
It is the smooth groupoid that will replace $TM$ in our tangent groupoid proof of the index theorem for contact manifolds.

As we saw, on a contact manifold all osculating groups $G_m$ are isomorphic to the Heisenberg group.
The converse is also true: if all osculating groups are Heisenberg groups, the manifold is a contact manifold (see [Er]).

If we identify $\UU(\Lg_m)$ with the algebra of right invariant differential operators on the Lie group $G_m$, then the principal part of $P$ at $m$ can be identified with a 
{\em homogeneous, right-invariant} differential operator $P_m$ on the osculating group $G_m$.
This operator is called the {\em model operator} of $P$ at $m$.
The principal symbol $\sigma^d_H(P)$ is thus identified with a smooth family of invariant differential operators on the fibers of $T_HM$.
Model operators on the Heisenberg group play the role in subelliptic theory on contact manifolds that constant coefficient operators on $\RR^n$ play in elliptic theory.
In the next section we will see that the main tool for the study of such operators is the same in both cases, namely, harmonic analysis.

\section{Rockland operators}

The basic idea of harmonic analysis on non-abelian groups is to replace Fourier theory with representation theory.
The appropriate generalization of an elliptic constant coefficient operator on $\RR^n$ in this context is a {\em Rockland operator}.

\begin{definition}
A {\em Rockland operator} on a graded nilpotent group $G$ is a differential operator $P$ that
is right-invariant, homogeneous, and has the property that $d\pi(P)$ is injective on the space of smooth vectors $\Sch_\pi$ for every irreducible unitary representation $\pi$ of $G$, except for the trivial representation.
\end{definition}
We recall some basic facts (see, for example, [Ta], [CG]).
A unitary representation $\pi$ of a Lie group $G$ on a Hilbert space $\HH_\pi$ induces a representation $d\pi$ of the Lie algebra $\Lg$ on the space of smooth vectors $\Sch_\pi\subseteq \HH_\pi$.
A {\em smooth vector} $v\in \HH_\pi$ is a vector for which the map $g\mapsto \pi(g)v$ is a $C^\infty$ function.
For $X\in \Lg$, the skew-Hermitian (generally unbounded) operator $d\pi(X)$ is defined by
\[ d\pi(X)v = \left. \frac{d}{dt}\right|_{t=0} \, \exp(tX) v ,\]
for $v\in \Sch_\pi$.
The representation $d\pi$ extends to the universal enveloping algebra $\UU(\Lg)$.
Thus, for every right invariant operator $P\in\UU(\Lg)$ we have an unbounded operator $d\pi(P)$ on $\HH_\pi$ with domain $\Sch_\pi$.

For the Heisenberg group $G$, the representation theory  is simple.
First of all, there is a one-parameter family of infinite dimensional irreducible representations $\pi_\tau$, $\tau\in \RR\setminus \{0\}$.
At the level of the Lie algebra $\Lg$ these representations are given by the formulas
\[ d\pi_\tau(X_j) = i\tau s_j \;,\; d\pi_\tau(Y_j) = \frac{\partial}{\partial s_j} \;,\; d\pi_\tau(T) = i\tau ,\]
where $s=(s_1,\ldots,s_n)\in \RR^n$, and $\HH_\pi = L^2(\RR^n)$.  
The smooth vectors $\Sch_\pi$ are the Schwartz class functions on $\RR^n$.

Additionally, there are the scalar representations of $G$, which are the characters of the abelian quotient $G/N\cong \RR^{2n}$, where the center $N$ is the set of elements $\{(0,0,t)\}\subseteq G$.
For all scalar representations we have $\pi(T)=0$.

\vskip 6pt
\noindent {\bf Examples.}
On an abelian group, a Rockland operator is just a homogeneous, constant coefficient {\em elliptic} operator. 
On the {\em graded} abelian group $\RR^{n+1}=\RR^n\times \RR$ the heat operator $P = \Delta + \partial/\partial t$
is a Rockland operator. 
Recall that the heat operator, even though it is not elliptic, is still `hypoellipitic'. (An operator $P$ is called hypoelliptic if, whenever $Pu=v$ for two distributions $u,v$, and $v$ is smooth in some open set, then $u$ is smooth in that same open set.)
This is true for Rockland operators in general.

As an instructive exercise the reader should verify that the operator
\[ P = \sum_{i=1}^n (-X_i^2-Y_i^2) + i\alpha T ,\]
on the Heisenberg group is a Rockland operator if and only if the complex constant $\alpha$ is not in the exceptional set $\{\ldots, -n-4,-n-2,-n,n,n+2,n+4,\ldots\}$.
In particular, the {\em sublaplacian}
\[ \Delta =  \sum_{i=1}^n (-X_i^2-Y_i^2) \]
is a Rockland operator.

This example illustrates quite dramatically that the `lower order term' of $P$ (in the ordinary sense) in the direction transversal to $H$ plays a crucial role in deciding whether or not the operator is regular.
This is the key observation underlying the Heisenberg calculus.
By assigning order $2$ to the vector field $T$, the term $i\alpha T$ becomes part of the principal symbol of $P$.
\vskip 6pt

Rockland operators satisfy {\em `a priori' estimates}, analogous to the well-known estimates for elliptic constant coefficient operators on $\RR^n$.
In full generality, this result is due to Helffer and Nourrigat ([HN]).

\begin{proposition}\label{hypind15:proposition:est}
If $P$ is a Rockland operator of degree $d$ on a graded group $G$,
and $A$ is an invariant differential operator of Heisenberg order $\le d$, 
then there exists a constant $C>0$ such that
\[ \|Au\| \le C \left(\| Pu\| + \| u\| \right),\]
for all $u \in \Sch(G)$. The norms in the inequality are $L^2(G)$ norms.
\end{proposition}
In our thesis [Er] we showed that the restrictions on the degree of the operator $P$ found in [HN] can be removed. We also verified that the theorem holds for {\em vector bundle} operators ([HN] only deals with scalar operators).

\section{Sobolev theory}

The a priori estimates for Rockland operators suggests an adaptation of {\em Sobolev theory} on the Heisenberg group that is compatible with the grading.
In this section we define the appropriate norms, and derive a crucial lemma.

Let $Y_0, Y_1,\cdots,Y_{2n}$ denote a basis for the Lie algebra of the Heisenberg group $\Lg$,
with $Y_0$ of degree $2$, and $Y_i$ of degree $1$ for $i=1,\ldots,2n$.
Then the {\em weigthed Sobolev space} $W^k = W^k(G)$ is the completion of $C^\infty_c(G)$
with respect to the norm,
\[ \|u\|^2_{W^k} = \sum_{|\alpha|\le k} \|Y^\alpha u\|^2 .\]
As usual, these weighted Sobolev spaces are Hilbert spaces.
The $W^k$ norms are, up to equivalence, independent of the choice of basis $Y_i$.

The a priori estimates for a Rockland operator $P$ can now be written in the usual form,
\[ \|u\|_{W^d} \le C(\|Pu\| + \|u\|) .\]
Because the expression $\|Pu\| + \|u\|$ is equivalent to the {\em graph norm} of $P$,
it follows that the domain of the closure of $P$ is the weighted Sobolev space $W^d(G)$.

To prove one of the key lemmas in this weighted Sobolev theory, 
we need an important result by Nelson and Stinespring ([NS], Theorem 2.2). 
The theorem in [NS] deals with elliptic operators, but its proof is valid for Rockland operators, because it relies only on the fact that the operator $P^*P+1$ is hypoelliptic.

\begin{theorem}
A formally self-adjoint Rockland operator $P$ on a graded group $G$ is essentially self-adjoint.
Moreover, if $\pi$ is any unitary representation of $G$,
then the closure of the operator $d\pi(P)$, with domain $\Sch_\pi$, is self-adjoint.
\end{theorem}
We now prove the analog of a well-known lemma in ordinary Sobolev theory.
It will be used several times in what follows.

\begin{lemma}\label{hypindex:lemma:sobolevplus}
Let $\Delta$ denote the sublaplacian on the Heisenberg group. 
The weighted Sobolev norm $\|u\|_{W^k}$ is equivalent to 
the norm $\|u\|_k = \|(\Delta+1)^{k/2} u\|$. 
\end{lemma}
{\bf Proof.}
If $k=2m$ is even, clearly $\|(\Delta + 1)^{m} u\| \le C\|u\|_{W^{2m}}$,
while the a priori estimates for $\Delta^m$ give us, 
\[ \|u\|_{W^{2m}} \le C\|(\Delta^{2m} + 1)^{1/2} u\| .\]
By spectral theory, $\|(\Delta^{2m} + 1)^{1/2} u\| \le \|(\Delta + 1)^{m} u\|$.
This proves the lemma for even $k$.
 
To deal with the odd case $k=2m+1$, simply observe that,
\[
 \|u\|^2_{W^{2m+1}} = \sum_{j=1}^n (\|X_j u\|^2_{W^{2m}}+\|Y_j u\|^2_{W^{2m}}) + \|u\|^2_{W^{2m}} 
                    = \ang{(\Delta+1) u,u}_{W^{2m}} 
                    = \|(\Delta+1)^{1/2}u\|^2_{W^{2m}} .
\]

\hfill{} $\Box$

As a corollary, we obtain a useful lemma.

\begin{corollary}\label{hypind15:lemma:epsilonC2}
Let $k\ge 1$. For every $\varepsilon>0$ there exists a $C>0$ such that
\[ \|u\|_{W^{k-1}} \le \varepsilon \|u\|_{W^k} + C  \|u\| .\] 
\end{corollary}
{\bf Proof.}
By spectral theory,
\[ \|(\Delta+1)^{(k-1)/2} u\| \le \epsilon \|(\Delta+1)^{k/2} u\| + C\|u\| .\]

\hfill{} $\Box$

\section{Subelliptic operators}

Let $\{U_j,\psi_j\}$ be an atlas for $M$, with Darboux coordinates $\psi_j\colon U_j \to G$, where $G$ is the Heisenberg group. 
If $T,X_1,\ldots,X_n,Y_1,\ldots,Y_n$ is the standard basis for the Lie algebra of right invariant vector fields on $G$, then
let $X^{(j)}_0,\cdots, X^{(j)}_{2n}$ denote the pullback of this basis to $U_j$. 
Let $\{\varphi_j\}$ be a partition of unity subordinate to $\{U_j\}$.
For a positive integer $k$, we define the {\em weighted Sobolev space} $W^k = W^k(M,H)$ on the contact manifold $M$
as the completion of $C^\infty(M)$ with respect to the norm,
\[ \|u\|^2_{W^k} = \sum_j \sum_{|\alpha|\le k} \|\varphi_j X_{(j)}^\alpha u\|^2 .\]
Again, these weighted Sobolev spaces are Hilbert spaces,
and the Sobolev norms are, up to equivalence,
independent of the choice of atlas $\{U_j\}$, or partition $\varphi_j$.
Obviously, a differential operator $A$ of Heisenberg order $d$ is continuous as an operator $W^{d+k}\to W^k$.

From Corollary \ref{hypind15:lemma:epsilonC2} one easily derives the analog for the Sobolev theory on contact manifolds.

\begin{lemma}\label{hypind15:lemma:sobolev}
Let $(M,H)$ be a compact contact manifold.
For every integer $k$, and every $\varepsilon>0$, there exists $C>0$ such that,
\[ \|u\|_{W^{k-1}} \le \varepsilon \|u\|_{W^k} + C\|u\| ,\]
for any $u\in C^\infty(M)$.
\end{lemma}
With this preparation, we can derive the main result about Fredholm operators in the Heisenberg calculus.
The reader should notice that the proof below would be {\em exactly} the same if we were deriving the a priori estimates for an elliptic operator!

\begin{theorem}\label{hypind15:thm:maxhypo}
Let $(M,H)$ be a compact contact manifold,
and let $P$ be a differential operator $P$ on $M$ for which all model operators $P_m$ are Rockland operators.
Then $P$ satisfies the a priori estimates
\[ \|u\|_{W^d} \le C(\|Pu\|+\|u\|) ,\]
for $u\in C^\infty(M)$.
\end{theorem}

\noindent {\bf Proof.}
By approximating the differential operator $P$ with its model operator $P_m$
in a neighborhood of each $m\in M$, we first show that
{\it for each $m\in M$ there is a neighborhood $V$ of $m$, and a constant $C_V$, such that,
\[ \|u\|_{W^d} \le C_V (\|Pu\| + \|u\|) ,\]
for functions $u$ with support in $V$.}

Let $U$ be a neigborhood of $m$, 
equiped with an isomorphism of Heisenberg structures $\psi\colon U\to G$,
where $G$ is the Heisenberg group. 
The Rockland operator $P_m$ now acts on functions $u\in C_c^\infty(U)$, and we have,
\[   \|u\|_{W^d} \le C(\|P_m u\| + \|u\|) .\]
All we need to do is compare $\|Pu\|$ and $\|P_m u\|$.

Choosing coordinates $(x_0,x_1,\cdots,x_{2n})$ in $U$, with $x=0$ at $m$, we have,
\[ P = P_m + \sum x_j Q_j + S,\]
where $Q_j$ are order $d$ differential operators, and $S$ is of order $d-1$. 
If $u$ is supported in a ball of small radius $|x|<\epsilon$, we get,
\[ \|P_mu\| - \|Pu\| \le \|(P_m-P)u\|   \le \epsilon C \| u \|_{W^d} + C \|u\|_{W^{d-1}} ,\]
with $C$ independent of $\epsilon$.  
Using Lemma \ref{hypind15:lemma:sobolev}, we obtain the desired estimate for smooth functions supported in an $\epsilon$ neighborhood of $m$, for $\epsilon$ sufficiently small.

To get the global result, choose a finite open cover $\{V_j\}$ 
such that the above estimates hold locally for each $V_j$.
Choose a smooth partition of unity $\{\varphi_j\}$ subordinate to $\{V_j\}$.
Then,
\begin{align*}
   \|Au\| \le \sum \|A\varphi_j u\| 
    & \le \sum C_j(\|P\varphi_j u\| + \|\varphi_j u\|) \\
    & \le \sum C_j(\|\varphi_j P u\| + \|\,[P, \varphi_j]\, u\| + \|\varphi_j u\|) \\ 
    & \le C(\|Pu\| + \|u\|) + \sum C_j \|\,[P, \varphi_j]\, u\|. 
\end{align*}
Again, Lemma \ref{hypind15:lemma:sobolev} finishes the proof.

\hfill $\Box$

The estimates in Theorem \ref{hypind15:thm:maxhypo} imply the weaker, subelliptic estimates
\[ \|u\|_{d/2} \le C(\|Pu\| + \|u\|) ,\]
where $\|\cdot\|_{d/2}$ denotes the {\em ordinary} Sobolev norm.
The sharper estimates in Theorem \ref{hypind15:thm:maxhypo} are sometimes refered to as {\em maximally hypoelliptic} estimates.
We need these sharper estimates to conclude that $P$ is Fredholm.

Observe that $\|Pu\|+\|u\|$ is equivalent to the norm on the graph of $P$.
The a priori estimates imply that the closure $\bar{P}$ of $P$ has domain $W^d$,
while the estimates extend by continuity to all $u\in W^d$. 
We need a regularity result.

\begin{proposition}\label{prop:moll}
Let $P$ be a differential operator with Rockland model operators.

If $u,v\in L^2$, and $Pu=v$ weakly (i.e., in the sense of distribution theory),
then $u$ is in the domain of the closure of $P$, and $\bar{P}u=v$.
\end{proposition}
The proof is a standard argument in elliptic theory, involving Friedrich's mollifiers.

\begin{corollary}
The symmetric operator 
\[ D = \left(\begin{array}{cc} 0 & P \\ P^t & 0 \end{array}\right), \]
acting on sections in $E\oplus F$, is essentially self-adjoint (i.e., its closure is self-adjoint).
\end{corollary}
{\bf Proof.}
It is a fact that
the formal adjoint of a  Rockland operator is also a Rockland operator 
Thus, the formal adjoint $P^t$ of $P$ is maximally hypoelliptic as well.
Therefore, by proposition \ref{prop:moll}, the closure of $P^t$ is the same as the adjoint $P^*$
(since $P^tu=v$ weakly and $P^*u=v$ mean the same thing for $u,v\in L^2$). 

\QED

\begin{theorem}
Let $P$ be a differential operator with Rockland model operators.
Then the closure of $P$ is Fredholm.
\end{theorem}
{\bf Proof.}
We can write the a priori estimates for the self-adjoint operator $\bar{D}$ as
\[ \|u\|^2_{W^d} \le C(\ang{\bar{D}u,\bar{D}u}+\ang{u,u}) = C\|(\bar{D}^2+1)^{1/2}u\|^2.\]
Thus, $(\bar{D}^2+1)^{-1/2}$ is a bounded operator $L^2\to W^d$, and therefore compact as an operator on $L^2$ (by the Rellich lemma).
It follows that $\bar{D}^2+1$ has discrete spectrum with finite dimensional eigenspaces.
In particular, the kernel of $\bar{D}$ (and therefore of $\bar{P}$ and $P^*$) is finite dimensional.
 
\QED

The index problem that presents itself is to determine the Fredholm index of a subelliptic operator on a compact contact manifold from its Heisenberg symbol, which is the family of Rockland model operators.
In the remainder of this paper we formulate a solution to this problem.

\include{cntinxI2} 

\chapter{An index theorem for subelliptic operators.}\label{chapter:index}

\section{The parabolic tangent groupoid of a contact manifold}\label{section:THM}

Our proof of the index theorem for `elliptic' operators in the Heisenberg calculus is an adaptation of Connes' tangent groupoid proof of the Atiyah-singer index theorem [Co].
Our first task here is to construct a groupoid associated to a contact manifold $(M,H)$ that plays the role of tangent groupoid.
We call our groupoid the {\em parabolic tangent groupoid} of $M$, and denote it $\THM$. The basic idea is pretty simple, but the technical realization is rather tricky.

Algebraically, the parabolic tangent groupoid of $M$ is the disjoint union of two smooth groupoids
\[ \THM = T_HM \cup M\times M\times (0,1] .\]
Here $M\times M\times (0,1]$ is a parametrized family of pair groupoids $M\times M$,
while $T_HM$ denotes the smooth groupoid that is (algebraically) the disjoint union of the osculating groups
\[ T_HM = \bigcup_{m\in M} G_m .\]
The fibers $G_m=H_m \oplus N_m$, with $N=TM/H$, are graded nilpotent groups (see section \ref{section:calc}).

We  make $\THM$ into a smooth groupoid by glueing $T_HM$ to the family of pair groupoids, in much the same way that Connes glued the tangent bundle $TM$ to this family. 
Connes achieved this by `blowing up the diagonal' in $M\times M$ by a factor $s^{-1}$ as the parameter $s\in (0,1]$ approaches $0$. In the limit, a pair of elements in $M\times M$ converges to a tangent vector.

To construct the parabolic tangent groupoid, we `blow up the diagonal' in $M\times M$, not by a simple factor $s^{-1}$, but by means of the dilations $\delta_s$ that are naturally associated with the graded groups $G_m$,
\[ \delta_s\;\colon (h,n) = (s^{-1}h, s^{-2}n), \]
for $(h,n)\in H_m\oplus N_m = G_m$.
This is compatible with the Heisenberg order of vector fields, since a vector field in the transversal direction (whose principal symbol is a section in the line bundle $N=TM/H$) defines a second order operator.
Correspondingly, it is blown up by a factor $s^{-2}$.

The `blowing up' must take place in a coordinate system on $M\times M$.
For the resulting procedure to yield a well-defined smooth structure on $\THM$ we must restrict ourselves to a carefully selected type of coordinates. We have considered the issue of the correct choice of coordinates in detail in our Ph.D. thesis [Er].
In the present paper we will simply present one way of effecting the dilation of $M\times M$.

Let $U\subseteq M$ be an open set endowed with Darboux coordinates,
i.e., $U$  is identified (as a contact manifold) with an open subset in the Heisenberg group $G$ (with its canonical contact structure).
Then the sub-groupoid $\THU\subseteq \THM$, which is the union
\[ \THU = T_HU \cup U\times U\times (0,1] \]
can be identified with a subset of the groupoid
\[ \THG = T_HG \cup G\times G\times (0,1] .\]
We will define a smooth structure on $\THG$ which is invariant under contact diffeomorphisms $G\to G$. Thus, $\THU$ inherits this structure, and an open cover of $M$ by sets endowed with Darboux coordinates will provide a smooth structure on $\THM$.

\section{The parabolic tangent groupoid of the Heisenberg group}\label{section:THG}

We now turn to the construction of the parabolic tangent groupoid $\THG$ of the Heisenberg group $G$. It can be elegantly described as a {\em continuous transformation groupoid}.

Let $\alpha$ be the action of the Heisenberg group $G$ on the space $B=G\times [0,1]$ given by
\[ \alpha(g)(p,s) = (\delta_s(g)p,s) ,\; g, p \in G,\; s\in [0,1] .\]
Here $\delta_s$ denotes the dilation of $G$ that corresponds to the grading, 
\[ \delta_s(x,y,t) = (sx,sy,s^2t) ,\]
for $g=(x,y,t) \in \RR^{2n+1}=G$.
Recall that the transformation groupoid $B\rtimes_\alpha G$ can be represented as the set of triples
$(b',g,b)\in B\times G\times B$ with  $b' = \alpha(g)b$. Triples compose as $(b'',g',b')(b',g,b) = (b'',g'g,b)$.
The smooth structure on $B\rtimes_\alpha G$ is obtained by the identification 
\[ B\rtimes_\alpha G \to G\times B \;\colon\; (b',g,b)\mapsto (g,b) .\] 
Algebraically, the groupoid $B\rtimes_\alpha G$ decomposes as a union of groupoids
\[ B\rtimes_\alpha G \cong \bigcup_{p\in G}G_p \;\cup\; G\times G\times (0,1] .\]
Here $G_p=\{((p,0),g,(p,0))\}$ is just a copy of the group $G$, one for each $p\in G$.
For $s>0$, we can identify the triple $((p',s),g,(p,s))$ with the pair $(p',p)$ in the {\em pair groupoid} $G\times G$ (one copy for each $s\in (0,1]$).

We now identify $G=G_p$ at $s=0$ with the osculating group $T_pG = H_p\oplus N_p$ of the contact manifold $G$ at the point $p$, by means of the right invariant trivialization $T_pG = \Lg = T_0G$ composed with the exponential map $\Lg\to G$. 
We obtain an isomorphism of groupoids
\[ B\rtimes_\alpha G \cong T_HG \cup G\times G\times (0,1] .\]
Thus, the family of pair groupoids $G\times G\times (0,1]$ is glued to the bundle of osculating groups $T_HG$ to form a smooth groupoid with boundary.
Observe that a path $(a_s,b_s,s)$ in $G\times G\times (0,1]$ corresponds to the triple
$((a_s,s), \delta_s^{-1}(a_sb_s^{-1}), (b_s,s))$ in $B\rtimes_\alpha G$. 
Therefore it converges to the element $g\in G_p$ in the osculating group at $p$ (corresponding to the triple $((p,0),g,(p,0))$) if  
\begin{align*}
 &\lim a_s = \lim b_s =  p , \\
 &\lim \delta_s^{-1}(a_sb_s^{-1}) =  g, 
\end{align*}
as $s\to 0$.
If we replace $G$ with the abelian group $\RR^n$, 
and the dilation $\delta_s$ with scalar multiplication,
then this convergence of $(a_s,b_s,s)\to g$ becomes simply
\[ \lim \frac{a_s-b_s}{s} = g ,\]
which describes precisely the topology of the tangent groupoid of $\RR^n$ (see [Co]).

\begin{definition}
The {\em parabolic tangent groupoid} of the Heisenberg group $G$ is the groupoid $\THG=T_HG\cup G\times G\times (0,1]$ with the smooth structure inherited from its identification with the transformation groupoid $(G\times [0,1])\rtimes_\alpha G$.
\end{definition}
To justify the identification of $G_p=\{((p,0),g,(p,0))\}$ with the osculating group at $p$,
we consider the effect of a contact transformation $\phi\colon G\to G$,
by which we mean a diffeomorphism that preserves the hyperplane bundle $H\subseteq TG$.
For a contact transformation, the derivative of $\phi$ induces a map  
\[ D\phi \;\colon\; H\oplus N \to H\oplus N \]
which is an automorphism of the bundle of osculating groups $T_HG$. 
At $s\in (0,1]$ we have the obvious automorphism of the pair groupoid,
\[ \phi\times\phi \;\colon\; G\times G\to G\times G .\]
Hence, a contact transformation $\phi$ naturally induces an automorphism of the groupoid  $\THG = T_HG \cup G\times G\times (0,1]$.

\begin{proposition}\label{prop:THG}
The automorphism of $\THG$
induced by a  contact transformation of $G$ is a diffeomorphism.
\end{proposition}
The proof is a straightforward argument involving Taylor expansions.
We leave the details to the reader.
The only non-trivial ingredient is the following basic fact about Taylor expansions of contact transformations.

\begin{lemma}
Let $\phi\colon G\to G$ be a contact transformation, such that $\phi(0)=0$ and $D\phi(0)=1$.
Denote coordinates on $G$ as $x=(x_0,x_1,\ldots,x_{2n})\in \RR^{2n+1}$, and write $y=\phi(x)$.
Then
\[ \frac{\partial^2y_0}{\partial x_i\partial x_j} = 0 ,\]
for $i,j=1,\ldots,2n$.
In other words, for any $v\in H_0\subseteq T_0G$ we have
$D^2\phi(v,v)\in H_0$.
\end{lemma}
{\bf Remark.} A second order Taylor expansion gives the corollary 
\[ \delta_s^{-1}\phi \,\delta_sg = g + \OO(s) ,\]
for any $g\in G$. 
This is the  non-trivial ingredient in the proof of Proposition \ref{prop:THG}.   

\vskip 6pt
\noindent {\bf Proof.}
Let $X_0,X_1,\ldots,X_{2n}$ denote the standard basis of $\Lg$, identified with right invariant vector fields on $G$, such that $X_i(0) = \partial/\partial x_i$ at $0\in G$.
We write $Y_i = T\phi(X_i)$.
By assumption, $Y_1,\ldots, Y_{2n}$ are sections in $H$, and $Y_i(0) = \partial/\partial x_i$.

Now fix $v\in \RR^{2n+1}$, and let $c(t)$ be the integral curve on $G$ with $c(0)=0$ and $c'(t) = \sum v_i\tilde{X}_i$. Assume that $v_0=0$.
We find
\[ c''(0) = \sum v_iv_j \frac{\partial Y_i}{\partial x_j}(0) .\]
For $i=1,\ldots, 2n$, we have $Y_i-X_i\in \sm(H)$ and $Y_i(0)-X_i(0)=0$, Therefore $\partial_j (X_i - Y_i)(0) \in \sm(H)$ for $i,j=1,\ldots, 2n$.
In other words, the $\partial/\partial x_0$ component $\partial_j Y_i^N(0)$ of $\partial_j Y_i(0)$ is equal to $\partial_j X_i^N(0)$, which denotes the $\partial/\partial x_0$ component of $\partial_j X_i(0)$, for $i,j=1,\ldots,2n$.
Inspecting the explicit expressions for the invariant vector fields $X_i$ on $G$ we see that
\[  \partial_j X_i^N + \partial_i X_j^N = 0 ,\]
for $i,j=1,\ldots,2n$.
It follows that $c''(0)^N=0$, or
\[ c''(0) \in H_0. \]
Because, in coordinates, the exponential map $\Lg\to G$ is the identity,
the curve $a(t)=tv$ is an integral curve of the vector field $\sum v_iX_i$.
Therefore, $c(t) = \phi(tv)$, and we find
\[ c''(0) = D^2\phi(v,v).\]
This proves the proposition.

\QED


\section{The index map induced by the parabolic tangent groupoid}

Following the procedure in [Co], observe that the restriction of a continuous, compactly supported function on the parabolic tangent groupoid $\THM$ to the $s=0$ boundary $T_HM$ induces a $\ast$-homomorphism $C^*(\THM)\to C^*(T_HM)$ with a contractible kernel.
Therefore the induced map in $K$-theory
\[ K_0(C^*(\THM)) \to K_0(C^*(T_HM)) \]
is an isomorphism.
Restriction to $s=1$ on the other hand gives a map
\[ K_0(C^*(\THM)) \to K_0(C^*(M\times M)) \cong K_0(\KK(L^2M)) \cong \ZZ  ,\]
and therefore a map
\[ {\rm Ind_H} \;\colon\; K_0(C^*(T_HM)) \to \ZZ .\]
Connes observes that this map (in his case, $K_0(C^*(TM))=K^0(T^*M)\to \ZZ$) corresponds to the analytic index in $K$-theory for elliptic operators,
and points to a proof of the Atiyah-Singer index theorem based on this observation.
Details of such a proof are worked out, to some extent, in [Hi].

We mimic this proof in the context of the Heisenberg calculus.
But first we must construct a class in $K_0(C^*(T_HM))$ representing the Heisenberg symbol $\sigma_H(P)$ of a subelliptic operator $P$ (i.e., the family of Rockland operators $P_m$).
This requires some new ideas.
Once this is established we will prove that ${\rm Ind}_H$ is, in fact, equal to the analytic index map, i.e., 
\[ {\rm Index}\, P = {\rm Ind_H}([\sigma_H(P)]) .\]

\section{A $K$-theory class for the Heisenberg symbol}\label{section:Ksymbol}

The principal part of a subelliptic operator $P$ in the Heisenberg calculus is represented by the smooth family of Rockland operators $P_m$ on the osculating groups $G_m$.
In this section we define a class $[\sigma_H(P)]\in K_0(C^*(T_HM))$ that is the analog of the symbol class $[\sigma(P)]\in K^0(T^*M)$ of an elliptic operator.

Let $D_m$ be the self-adjoint operator,
\[ D_m = \left(\begin{array}{cc} 0 & -iP_m \\ iP_m^* & 0 \end{array}\right) ,\]
and let $u_m$ be its Cayley transform,
\[ u_m = (D_m+i)(D_m-i)^{-1} .\]
Note that $D_m$ is just the model operator of a self-adjoint operator $D$ associated to $P$.
All operators $D_m$ are Rockland operators.

If $P$ acts on sections in bundles $E$ and $F$, then $D_m$ acts on sections in the trivial bundle $\pi^*E\oplus\pi^*F$ over $G_m$. (Here $\pi$ is the base point map $\pi\colon T_HM \to M$.)
Let $\epsilon$ denote the grading operator associated to the decomposition $\pi^*E\oplus \pi^*F$,
\[ \epsilon = \left(\begin{array}{cc} -1 & 0 \\ 0 & 1 \end{array}\right). \]
We have $\epsilon^2 = 1$, and $\epsilon D_m = - D_m\epsilon$. 
It follows that $\epsilon u_m = u_m^* \epsilon$, and then $(\epsilon u_m)^2 = 1$. 
This last equality implies that $\frac{1}{2}(\epsilon u_m + 1)$ is a projection.

\begin{definition}
Let $P$ be a subelliptic differential operator on a compact contact manifold $(M,H)$,
with Rockland model operators.
Then the symbol class in $K$-theory associated to $P$ is defined as
\[  [\sigma_H(P)] = [\frac{1}{2}(\epsilon u+1)] - [\frac{1}{2}(\epsilon +1)] \in K_0(C^*(T_HM)) .\]
\end{definition}
This expression must be treated with care.
The projection,
\[ e = \frac{1}{2}(\epsilon u+1) ,\]
denotes the {\em family} of projections $\{e_m,m\in M\}$, where
\[ e_m = \frac{1}{2}(\epsilon u_m+1) \] 
acts on the Hilbert space $L^2(G_m)\otimes (E_m\oplus F_m)$.
We have to show first of all that
\[ e_m \in C^*(G_m)^+ \otimes {\rm End}\,(E_m\oplus F_m) \cong \Mat_{2N}(C^*(G_m)^+), \]
where $2N$ denotes the fiber dimension of $E\oplus F$.
Since $u_m = 1+2i(D_m-i)^{-1}$, it suffices to prove that the resolvent of the Rockland operator $D_m$ satisfies
\[ (D_m-i)^{-1} \in \Mat_{2N}(C^*(G_m)).\]
Secondly, we must prove that the family $\{e_m\}$ defines a continuous section,
\[ e \in C^*(T_HM)^+\otimes_M\, C({\rm End}(E\oplus F)) .\]
Here $C^*(T_HM)\otimes_M\, C({\rm End}(E\oplus F))$ denotes the $C^*$-algebra of continuous sections in the field $\{C^*(G_m)^+ \otimes {\rm End}\,(E_m\oplus F_m)\}$ (Morita equivalent to $C^*(T_HM)^+$).
By choosing a bundle $K$ such that $E\oplus F\oplus K$ is trivial, we obtain,
\[ C^*(T_HM)^+\otimes_M\, C({\rm End}(E\oplus F)) \subseteq \Mat_k(C^*(T_HM)^+) .\]
Again, it suffices to show that the family of resolvents $\{(D_m-i)^{-1}\}$ defines a continuous section in the field $\{\Mat_{2N}(C^*(G_m))\}$ (here we can simply choose a {\em local} trivialization of $E\oplus F$).

We now verify these facts.

\begin{proposition} \label{hypind15:prop:main}
Let $G$ be a graded nilpotent group.
If we identify $C^*(G)$ with its image under the left regular representation on $L^2(G)$,
then the resolvent of a self-adjoint Rockland operator on $G$ is an element in $C^*(G)$.
\end{proposition} 
{\bf proof.}
A theorem of Folland and Stein (see [FS]) states that, for the {\em positive} Rockland operator $L^2$,
the distribution kernel $K$ of $\exp(-L^2)$ (a kind of `heat kernel'),
\[ e^{-L^2}\,\phi = K\ast \phi, \]
is a Schwartz-class function on the group $G$ ([FS], chapter 4.B).
Hence, certainly $K\in L^1(G)$, and therefore,
\[  e^{-L^2} \in C^*(G) .\]
Next, the distribution kernel of the operator $L e^{-L^2}$ is $(L\delta) \ast K = LK$.
The function $LK$ is of Schwartz class, because any left invariant homogeneous operator $L = \sum a_\alpha(x)\partial^\alpha$ has polynomial coefficients $a_\alpha(x)$. 
We conclude that also, 
\[ L e^{-L^2} \in C^*(G) .\]
Because the functions $e^{-x^2}$ and $x e^{-x^2}$ separate points on $\RR$, these two functions generate $C_0(\RR)$ as a $C^*$-algebra.
Hence $f(L)\in C^*(G)$ for any $f\in C_0(\RR)$.

\hfill $\Box$

The theorem of Folland and Stein quoted here also holds if $L$ is a matrix of operators.
(Their proof relies only on the existence of a priori estimates for a sufficiently large power of $L$, and these a priori estimates, 
and therefore Proposition \ref{hypind15:prop:main}, also holds for such operators.)
Thus we have established the fact that,
\[ (D_m-i)^{-1} \in \Mat_{2N}(C^*(G_m)) .\]
Our second requirement is that the family $\{(D_m-i)^{-1}\}$ defines a {\em continuous section}
in the field of $C^*$-algebras $\{\Mat_{2N}(C^*(G_m))\}$ over $M$.
On a contact manifold, the field $\{\Mat_{2N}(C^*(G_m))\}$ is locally trivial,
and we can think of the model operators $D_m$ simply as a family of Rockland operators on a fixed group $G$, with smoothly varying coefficients.
All we need to show is that the map $m\mapsto (D_m-i)^{-1}$ for such a family is norm continuous.
This is a straightforward exercise, which uses the basic equality 
\[ (D_{m_0} - i)^{-1} - (D_{m_1} - i)^{-1} =
   (D_{m_1} - i)^{-1}  (D_{m_1} - D_{m_0})  (D_{m_0} - i)^{-1} .\]
One considers the norms of the three factors on the right-hand side as operators (from left to right) from $L^2$ to $L^2$, from $W^d$ to $L^2$, and from $L^2$ to $W^d$ (by virtue of subellipticity of $D_{m_0}$). As $m_1\to m_0$, the norm of the second factor converges to zero, while the other two remain uniformly bounded.

This establishes that $[\sigma_H(P)]$ is a well-defined element in $K_0(C^*(T_HM))$.

\vskip 6pt
\noindent {\bf Remark.}
If we take $H=TM$, and an elliptic operator $P$ with constant coefficient operators $P_m$ on $T_mM$, the construction described in this section results in an element in  $K_0(C^*(TM))$.
This element corresponds to the usual element
$[\sigma(P)] \in K^0(T^*M)$
in topological $K$-theory.
The key fact is that the projection $e_m$ is identical to the {\em graph projection} of $P_m$.

\section{The Fredholm index as a $K$-theory class}

The construction described in the previous section can be applied to $P$ itself.
We start with the self-adjoint operator 
\[ D = \left(\begin{array}{cc} 0 & -i\bar{P} \\ iP^* & 0 \end{array}\right) \]
on the Hilbert space $\HH = L^2(M,E)\oplus L^2(M,F)$,
with Cayley transform 
\[ u = (D+i)(D-i)^{-1} .\]
(Here $\bar{P}$ denotes the closure of $P$.)
The direct sum decomposition of the Hilbert space $\HH$ is encoded in the grading operator $\epsilon$.
As before, we obtain a $K$-theory element
\[ [e_P] - [f_P] = [\frac{1}{2}(\epsilon u  + 1)] - [\frac{1}{2}(\epsilon +1)] \in K_0(\KK(\HH)).\] 
Because $ u = 1 + 2i(D-i)^{-1}$ we have
\[ e_P - f_P = \frac{1}{2}(\epsilon u  + 1) - \frac{1}{2}(\epsilon +1) = i\epsilon(D-i)^{-1}  ,\]
and so $e_P-f_P$ is compact if the resolvent $(D-i)^{-1}$ is compact, which holds when $D$ is subelliptic.

\begin{proposition}
With the usual identification $K_0(\KK)\cong \ZZ$, we have
\[ [e_P]-[e_P] = {\rm Index}\, P.\]
\end{proposition}
{\bf Proof.}
consider the family of projections $e_t = \frac{1}{2}(\epsilon u_t + 1)$,
with $u_t = (tD+i)(tD-i)^{-1}$.
By functional calculus, the family $u_t$ is norm continuous.
Observe that the spectral function $z\mapsto (z+i)(z-i)^{-1}$ maps $0$ to $-1$, and $\infty$ to $+1$.
Because $D$ has compact resolvent it has discrete spectrum.
Therefore, as $t\to \infty$, the family $u_t$ converges in norm to the operator,
\[ u_\infty = -[{\rm Ker}\,D] + (1-[{\rm Ker}\,D]) = 1 - 2[{\rm Ker}\,D].\]
It now follows from homotopy invariance of $K$-theory that
\[ [e_P] - [f_P] = [{\rm Ker}\,P] - [{\rm Ker}\, P^*] \in K_0(\KK(\HH)),\]
where $[{\rm Ker}\,P], [{\rm Ker}\,P^*]$ denote the projections of $\HH$ onto the kernels of $P$ and $P^*$.

\QED

A similar homotopy in the opposite direction, i.e., with $t\to 0$, will continuously deform $[e_P]-[f_P]$ into $[\sigma_H(P)]$.
This, in essence, is the proof of the index theorem. 
We will now make this precise.

\section{The analytic index map for subelliptic operators}\label{section:aind}

Let $P$ be a differential operator on $M$ of Heisenberg degree $d$.
We construct a right invariant differential operator $\PP$ on the groupoid $\THM$,
and an element in $K_0(C^*(\THM)$ associated to $\PP$.
For simplicity of exposition, we will describe the construction for the case of scalar operators $P$.
However, it is a trivial matter to adapt the construction in the case where $P$ is   vector bundle operator $P\;\colon\; \sm(E) \to \sm(F)$.

Recall that the source map of a smooth groupoid is a submersion, and the source fibers (the pre-image of points in the base space) are smooth manifolds.
Right multiplication by an arrow in the  groupoid induces a diffeomorphism between two source fibers.
A right invariant operator is a differential operator on the groupoid that differentiates in the direction of the source fibers, and is invariant with respect to right multiplication.
(see [Ma])

Every osculating group $G_m$ is a source fiber, and on $G_m$ we let 
\[ \PP_m=P_m .\]
This means that at the $s=0$ boundary in the parabolic tangent groupoid $\THM$,
the operator $\PP$ is just the smooth family of model operators $P_m$ on the bundle of osculating groups $T_HM$.
The source fibers of the pair groupoid $M\times M$ are just copies of $M$,
and right invariance simply means that we have the same operator in each fiber.
In a $s>0$ source fiber of $\THM$ we let 
\[ \PP_s=s^dP .\]
If we pick an open set $U\subseteq M$ with Darboux coordinates $U\to G$,
and identify $\THU$ as a subset of $\THG=B\rtimes_\alpha G$ , then $\PP_s$ is represented as the operator
\[ \PP_s = P_d + sP_{d-1} + s^2P_{d-2} + \cdots \]
on the $s$-fibers $G$ of the action groupoid $B\rtimes_\alpha G$ that defines the smooth structure of $\THG$ (see section \ref{section:THG}).
Here $P_k$ denotes the part of $P$ that is homogeneous of Heisenberg order $k$.
We see that $\PP$ has smooth coefficients on $\THM$.

Let $D$ be the self-adjoint operator associated to $P$, as before,
and let ${\mathbb D}$ be the right invariant operator on $\THM$ associated to $D$, in the same way that $\PP$ is associated to $P$.
The key technical result is the following theorem. We postpone its proof until later.

\begin{proposition}\label{hypind50:prop:bigD}
Let $(M,H)$ be a compact contact manifold,
and $D$ a self-adjoint differential operator on $M$ of Heisenberg order $d$ with Rockland model operators.
Then 
\[ ({\mathbb D} - i)^{-1}\in \Mat_k(C^*(\THM)) .\] 
\end{proposition}
{\bf Remark.}
The integer $k$ denotes the dimension of a trivial vector bundle over $M$ that contains $E\oplus F$ as a direct summand.
\vskip 6pt

Proposition \ref{hypind50:prop:bigD} immediately implies the following theorem.

\begin{theorem}\label{hypind50:thm:analindx}
Let $(M,H)$ be a compact contact manifold,
and $P$ a subelliptic operator on $M$ with Rockland model operators.
Then,
\[ {\rm Index}\,P = {\rm Ind}_H\,([\sigma_H(P)]) .\]
\end{theorem}
{\bf Proof.}
Proposition \ref{hypind50:prop:bigD} implies that there is a well-defined unitary,
\[ U = ({\mathbb D} + i)({\mathbb D} - i)^{-1} = 1-2i({\mathbb D} - i)^{-1} \in \Mat_k(C^*(\THM)^+) ,\]
and a projection,
\[ \frac{1}{2}(\epsilon U + 1)\in  \Mat_k(C^*(\THM)^+) .\]
We thus obtain a $K$-theoretic `index' for the invariant family $\PP$,
\[ [\PP] = [\frac{1}{2}(\epsilon U + 1)] - [\frac{1}{2}(\epsilon + 1)] \in K_0(C^*(\THM)) .\]
Its restrictions to $s=0$ is $[\sigma_H(P)] \in K_0(C^*(T_HM))$, and its restriction to $s=1$ is the element in $K_0(\KK(L^2(M)))$ that corresponds to the Fredholm index of $P$.

\hfill $\Box$

\section{The topological index for subelliptic operators}\label{section:tind}

In this section we will prove the main result of this paper.

\begin{theorem}\label{MAIN}
Let $(M,H)$ be a compact contact manifold,
and $P$ a subelliptic differential operator on $M$ with Rockland model operators.
Then
\[ {\rm Index}\,P = \int_{T^*M} \, {\rm Ch}(\Psi([\sigma_HP]))\,\wedge\, {\rm Td}(M) .\]
\end{theorem}
In other words, the analytic index of $P$ is computed by the topological index formula of Atiyah and Singer.
Here 
\[ \Psi\;\colon\; K_0(C^*(T_HM)) \stackrel{\cong}{\longrightarrow} K^0(T^*M) \]
denotes a well-known {\em canonical} isomorphism in $K$-theory that we will describe below.
\vskip 6pt

Consider the {\em adiabatic groupoid} of $\THM$ (see, for example, [Ni2]). 
Recall that we can think of $\THM$ as a family of groupoids over the unit interval $[0,1]$, where at $t=0$ we have the bundle of osculating groups $T_HM$, while at each $t>0$ we have a copy of the pair groupoid $M\times M$.
If we analyze the adiabatic groupoid $\THM^{ad}$, we see that it is the union of a family of groupoids $\GG_{(t,s)}$ parametrized by $(t,s)\in [0,1]^2$,
\begin{align*}
\GG_{(t,s)} &= M\times M\; {\rm for}\;t>0,s>0,\\
\GG_{(0,s)} &= T_HM ,\; {\rm for}\;s>0\\
\GG_{(t,0)} &= TM,\; {\rm for}\;t>0\\
\GG_{(0,0)} &= H\oplus N.
\end{align*}
Each groupoid $\GG_{(t,s)}$ has unit space $M$, and the unit space of $\THM^{ad}$ is the manifold with corners $M\times [0,1]^2$.  
In what follows it may be useful to refer to the following schematic picture,
which sums up how the groupoids $\GG_{(t,s)}$ are assembled along the horizontal and vertical edges in the square $[0,1]^2$.
\[ \xymatrix{ {H\oplus N} \ar@{.}[r]|-{\delta_t^{-1}}  \ar@{.}[d]_{s^{-1}} &  {TM} \ar@{.}[d]_{s^{-1}}^{\TM} \\
              {T_HM}      \ar@{.}[r]^-{\THM}_-{\delta_t^{-1}} & {M\times M}     }
\]
The groupoid $\THM^{ad}$ gives rise to a commutative diagram in $K$-theory, induced by restriction of functions on $\THM^{ad}$ to each of the four corners of the square $[0,1]^2$. We proceed step-by-step.

\begin{lemma}
Restriction of elements in $C^*(\THM^{ad})$ to the $(t,s)=(0,0)$ corner,
\[ e_{(0,0)} \;\colon\; C^*(\THM^{ad}) \to C^*(H\oplus N)\cong C_0(H^*\oplus N^*), \]
induces an isomorphism in K-theory.
\end{lemma} 
{\bf Proof.}
Let $\GG$ denote the groupoid that is the union of the $t=0$ and $s=0$ edges in $\THM^{ad}$.
The restriction map $C^*(\THM^{ad})\to \GG$ induces an isomorphism in $K$-theory,
because the kernel of this map is the contractible ideal $C_0((0,1]^2,\KK)$.
The kernel of the map that further restricts $\GG$ to the corner $(t,s)=(0,0)$ is again a contractible ideal. 

\QED 

Now let $e_{(t,1)}$ denote restriction to the edge $s=1$,
\[ e_{(t,1)}\;\colon\; C^*(\THM^{ad}) \to C^*(\THM) ,\]
and $e_{(1,s)}$  restriction to the edge $t=1$,
\[ e_{(1,s)}\;\colon\; C^*(\THM^{ad}) \to C^*(\TM) .\]
Further restriction to the corner $(t,s)=(1,1)$ gives two $\ast$-homomorphisms,
\begin{align*}
e_{(1,1)}&\;\colon\; C^*(\THM) \to C^*(M\times M) \\
e'_{(1,1)}&\;\colon\; C^*(\TM) \to C^*(M\times M).
\end{align*}
We obtain a commutative diagram,
\[ \xymatrix{  K_0(C^*(\THM^{ad})) \ar[d]^{e_{(t,1)}} \ar[r]^{e_{(1,s)}} & K_0(C^*(\TM)) \ar[d]^{e'_{(1,1)}} \\
               K_0(C^*(\THM)) \ar[r]_-{e_{(1,1)}}& K_0(C^*(M\times M)) }
\]
But $e_{(1,1)}$ is just our deformation index ${\rm Ind}_H$,
while $e'_{(1,1)}$ is the deformation index associated to the usual tangent groupoid.
As is known from the tangent groupoid proof of the Atiyah-Singer index theorem for elliptic operators, this map $e'_{(1,1)}$ is equal to the topological index ${\rm Ind_t}$ (see [Co], [Hi]).
Therefore, the diagram above can be read as follows,
\[ \xymatrix{  K^0(H^*\oplus N^*) \ar[d]^{e_{(t,1)}} \ar[r]^-{e_{(1,s)}} & K^0(T^*M) \ar[d]^{\rm Ind_t} \\
               K_0(C^*(T_HM)) \ar[r]_-{\rm Ind_H} & \ZZ }
\]
The remaining maps $e_{(1,s)}$ and $e_{(t,1)}$ are natural isomorphisms.
For $e_{(1,s)}$ this is trivial.
But for $e_{(t,1)}$ this result is far from trivial.

\begin{proposition}\label{hypind50:prop:et1}
The map
\[ e_{(t,1)}\;\colon\; K^0(H^*\oplus N^*) \to K_0(C^*(T_HM)) \]
is an isomorphism in $K$-theory.
\end{proposition}
This is a special case of Lemma 3 in [Ni2].  
We sketch the proof, to show how it relies on the highly nontrivial {\em Connes-Thom Isomorphism}.
The osculating group $G_m$, being nilpotent, can be represented as a semi-direct product $G_m\cong G'_m \rtimes \RR$, with $G'_m$ again nilpotent.
We get a crossed product
\[ C^*(G_m) \cong C^*(G'_m)\rtimes \RR .\]
By the Connes--Thom isomorphism for crossed products with $\RR$ (see [Co1]), we have an isomorphism
\[ K_0(C^*(G_m)) \cong K_0(C^*(G')\otimes C_0(\RR^*)).\]
Inductively, we obtain {\em natural} isomorphisms
\[ K_0(C^*(G_m)) \cong K^0(\Lg_m^*) ,\]
where $\Lg_m^*$ is the dual space of the Lie algebra $\Lg_m = H_m\oplus N_m$ of the osculating group $G_m$.
It is  not hard to show that the map $e_{(t,1)}$ restricts to the Connes--Thom isomorphism for the osculating group $G_m$ in each of the fibers of $T_HM$. 
Once this is established, the proof of Proposition \ref{hypind50:prop:et1} is completed by a Mayer--Vietoris argument. (See [Ni2] for details.)
\vskip 6pt

What results is a natural isomorphism
\[ \Psi = e_{(1,s)}\circ e_{(t,1)}^{-1} \;\colon\; K_0(C^*(T_HM)) \stackrel{\cong}{\longrightarrow} K^0(T^*M) \]
between the $K$-theory group that receives our noncommutative Heisenberg symbol and the topological $K$-theory of $T^*M$, and our commutative diagram reduces to
\[ \xymatrix{   & K^0(T^*M)  \ar[d]^{\rm Ind_t} \\
               K_0(C^*(T_HM)) \ar@/^1pc/[ur]^{\Psi}  \ar[r]_-{\rm Ind_H} & \ZZ}
\]
In particular,
\[ {\rm Ind}_H([\sigma_HP]) = {\rm Ind_t}(\Psi([\sigma_HP])).\]
In combination with Theorem \ref{hypind50:thm:analindx}, this proves Theorem \ref{MAIN}.

\section{A technical result}
\label{hypind40:sect:C*THMfamilies}

Our aim in this section is to prove Proposition \ref{hypind50:prop:bigD},
which is the `hard nut' hidden inside the soft and elegant machinery of groupoids and commutative diagrams.

We start by considering the convolution $C^*$-algebra $C^*(\THG)$ of the parabolic tangent groupoid $\THG$ of the Heisenberg group, and derive a concrete characterization of the elements in this $C^*$-algebra.

To understand the structure of $C^*(\THG)$ it is useful to consider the various restriction maps associated to the stratification of the groupoid.
Restriction at $s=0$ of functions in the convolution algebra $C_c(\THG)$ determines a $*$-homomorphism to the convolution algebra $C_c(T_HG)$ which extends by continuity to a surjective $*$-homomorphism of the corresponding $C^*$-algebras, 
\[ \pi_0 \;\colon\; C^*(\THG) \to  C^*(T_HG) .\]
Further restriction to the individual osculating groups gives morphisms
\[ \pi_p \;\colon\; C^*(\THG) \to C^*(G_p) .\]
Likewise, for $s\in (0,1]$ we obtain surjective $*$-homomorphisms,
\[ \pi_s \;\colon\; C^*(\THG) \to C^*(G\times G) \cong \KK(L^2(GG) .\]
Thus, every element $Q\in C^*(\THG)$ gives rise to a family of elements
\[ \ang{Q_p,Q_s} = \ang{\pi_p(Q),\pi_s(Q)} ,\]
with 
\[ Q_p\in C^*(G) ,\; Q_s\in \KK(L^2(G)) .\] 
Our aim in this section is to determine precisely which families $\ang{Q_p,Q_s}$ of this type represent elements in $C^*(\THG)$. 

We derive two properties.
In what follows $f_s$ denotes the dilated function $f_s(x) = f(\delta_sx)$ on $G$.

\begin{proposition}\label{hypind40:prop:asymplocal}
If $Q=\ang{Q_p,Q_s}\in C^*(\THG)$, then for two continuous functions $\phi, \psi\in C(M)$ with disjoint supports we have
\[ \lim_{s\to 0} \, \| \phi_s \,Q_s\, \psi_s \| = 0.\]
Here $\phi_s(x) = \phi(\delta_sx)$.
\end{proposition}
{\bf Proof.}
This is trivial for kernels $k\in C_c(\THG)$
The result for arbitrary $Q\in C^*(\THM)$ follows by approximation.

\QED

The second property of elements $\ang{Q_p,Q_s}$ is less obvious. 
It describes how the family $Q_s$ converges to the family $Q_p$ as $s\to 0$.
It establishes the sense in which an operator in the Heisenberg calculus deforms to its Heisenberg symbol.

\begin{proposition}\label{hypind40:prop:asymp2}
If $Q=\ang{Q_p,Q_s}\in C^*(\THG)$,
then for every  $p\in G$ and every $\varepsilon > 0 $ there exists a neighborhood $V$ of $p$ such that, if $h$ denotes the characteristic function of $V$, then
\[ \limsup_{t\to 0}\; \|\, h_s Q_s h_s - h_s  Q_m h_s \|  < \varepsilon .\]
\end{proposition}
{\bf Proof.}
We identify $\THG \cong (G\times [0,1])\rtimes G = B\rtimes G$,
and consider a kernel $k\in C_c((G\times [0,1])\rtimes G)$.
The regular representations of $k$ on $L^2(\GG_{(p,s)}) =L^2(G)$
are given, in triple notation, by
\[ \pi_{(p,s)}(k) f(x) = \int_G k(((\delta_sx)p,s),xy^{-1},((\delta_sy)p,s)) \,f(y) dy .\]
Here we write $f(y) = f((\delta_sy)b,y,b)$.
At $s>0$ we may identify this operator with $Q_s = \pi_s(k)$, while at $s=0$ we obtain the convolution operator $Q_p = \pi_p(k)$.
If $a_s(x,y)$ denotes the Schartz kernel of $h_s k_s h_s$ and $b_s(x,y)$ that of $h_s k_m h_s$, we have
\begin{align*}
 a_s(x,y) & = h(\delta_sx)\; k(((\delta_sx)p,s),xy^{-1},((\delta_sy)p,s))\; h(\delta_sy), \\
 b_s(x,y) & = h(\delta_sx)\; k((p,0),xy^{-1},(p,0))\; h(\delta_sy) .
\end{align*}
By taking the support $V\subseteq G$ of $h$ sufficiently small,
and letting $s<r$ for sufficiently small $r>0$, 
we can obtain a uniform estimate
\[ | a_s(x,y) - b_s(x,y) | < \varepsilon .\]
To see this, 
consider the values $(x,y,s)$ for which $a_s(x,y)-b_s(x,y)\ne 0$.
We must have $\delta_sy\in V$, but also $xy^{-1} \in K$, where $K\subseteq G$ is the compact set 
\[ K = \{ g\in G \;|\; (b',g,(v,s))\in {\rm supp}(k), \; v\in V, s\in [0,r] \} .\]
This set is compact because $k$ has compact support. 
Then, because $k$ is continuous, the difference 
\[ | k(((\delta_sx)p,s),xy^{-1},((\delta_sy)p,s)) - k((p,0),xy^{-1},(p,0)) | \]
can be made arbitrarily small for all $\delta_sy\in V, s<r$,
by taking $V$ and $r>0$ small enough.
This estimate will be {\em uniform} in $xy^{-1}\in K$ (because $K$ is compact),
which gives the required estimate on $|a_s-b_s|$.

Let us denote $c_t(x,y)=a_t(x,y)-b_t(x,y)$.
The fact that $xy^{-1}\in K$ if $c_t(x,y)\ne 0$ implies uniform estimates 
\[ \int |c_t(z,y)| dy < {\rm Vol}(K)\varepsilon ,\; \int |c_t(x,z)|dx < {\rm Vol}(K)\varepsilon ,\]
that hold for each $z\in G$. Hence,
\[ \|h_sk_sh_t - h_tk_m h_t\|  < {\rm Vol}(K)\varepsilon ,\]
for all $s<r$. This gives the required estimate.

The result for general $\ang{Q_p,Q_s}\in C^*(\THG)$ follows by approximation.

\QED

The two properties we derived are sufficient to characterize elements in $C^*(\THG)$.

\begin{proposition}\label{hypind40:prop:asymp3}
A family of operators $\ang{Q_p,Q_s}$ with $Q_p\in C^*(G)$ and $Q_s\in \KK(L^2(G))$ represents an element in $C^*(\THG)$ if and only if it has all of the following properties:

(1) the family $Q_p$ is norm continuous;

(2) the family $Q_s$ is norm continuous and uniformly bounded;

(3) the family $Q_s$ satisfies the condition of Proposition \ref{hypind40:prop:asymplocal}; 

(4) the family $\ang{Q_p,Q_s}$ satisfies the condition of Proposition \ref{hypind40:prop:asymp2}.

\end{proposition}
{\bf Proof.}
We have proven necessity of the properties listed, and must now show that they are sufficient.
Let $\DD$ be the set of elements $\ang{Q_p,Q_s}$ that satisfy all four properties.
Because of properties (1) and (2), we can think of $\DD$ as a subset of the $C^*$-algebra
\[ C^*(T_HG) \oplus C_b((0,1],\KK)) .\]
One easily verifies that $\DD$ is a norm-closed $*$-subalgebra, and therefore a $C^*$-algebra.
For example, to verify property (3) for the product, let $\ang{Q_p,Q_s}$ and $\ang{R_p,R_s}$ be two families in $\DD$. 
Let $\phi_1,\phi_2\in C_c(G)$ be two functions with disjoint support.
Choose two other function $\psi_1,\psi_2\in C_c(G)$, such that $\psi_1$ and $\psi_2$ have disjoint supports as well, and such that $\psi_j(x)=1$ whenever $x\in {\rm supp}(\phi_j)$.
Then $(1-\psi_j)$ and $\phi_j$ have disjoint supports as well, and therefore
\begin{align*}
  \lim_{s\to 0} & \| \phi_1 Q_s - \phi_1 Q_s \psi_1 \| = \lim_{s\to 0} \| \phi_1 Q_s (1- \psi_1) \| = 0, \\
  \lim_{s\to 0} & \| R_s \phi_2 - \psi_2 R_s \phi_2 \| = \lim_{s\to 0} \| (1 - \psi_2) R_s \phi_2 \| = 0.
\end{align*}
Because both $Q_s$ and $R_s$ are uniformly bounded in norm, it follows that
\[ \lim_{s\to 0} \| \phi_1 Q_tR_t \phi_2 \| = \lim_{t\to 0}  \| \phi_1 Q_t\psi_1\psi_2R_t \phi_2 \| = 0.\]
This establishes property (3) for the product.
Property (4) is proven in a similar way, while properties (1) and (2) are trivial.

We see that $\DD$ is indeed a $C^*$-algebra,
and by the previous two propositions we know that $C^*(\THG)\subseteq \DD$.
To see that the two are isomorphic, consider the restriction map
\[ \DD\to C^*(T_HG) \;\colon\; \ang{Q_p,Q_s}\mapsto \{Q_p\} .\]
Since we have a short exact sequence
\[ 0\to C_0((0,1],\KK) \to C^*(\THG) \to C^*(T_HG) \to 0 ,\]
we only need to show that the $s=0$ restriction for $\DD$ has the same kernel.

Suppose therefore that $\ang{Q_p,Q_s}\in \DD$ and that $Q_p=0$ for all $p\in G$.
Choose $\varepsilon >0$.  
By property (4) there exists a neigborhood $V$ of each point $p\in G$ such that the characteristic function $h$ of $V$ satisfies
\[ \limsup_{s\to 0} \|h_sQ_sh_s\| < \epsilon .\] 
This clearly implies $\|Q_s\| < \epsilon$, which proves the claim.

\QED

\vskip 6pt
\noindent{\bf Proof of Proposition \ref{hypind50:prop:bigD}.}
For convenience, we denote $f(x) = (x-i)^{-1}$.
We must verify four properties for the family $\ang{f(D_p),f(s^dD)}$.
Property (1) was established in Proposition \ref{hypind15:prop:main}.
Property (2) follows by functional calculus.
Property (3) follows from Lemma \ref{hypind99:lemma:commutator} below.

We must verify property (4)
\[  \limsup_{s\to 0} \|h_s f({\mathbb D}_p) h_s - h_s f({\mathbb D}_s) h_s \| < \varepsilon \]
on the source fibers $G$ in $\THG$ ($U\subseteq M$, $U\to G$),
or, equivalently,
\[  \limsup_{s\to 0} \|h f(s^dD_m) h - h f(s^dD) h \| < \varepsilon \]
on the source fibers $M$ of $M\times M$.
The standard trick for comparison of two commutators gives us
\[ (s^dD_m - i)^{-1} - (s^dD-i)^{-1} = s^d\,(s^dD_m - i)^{-1}(D-D_m)(s^dD-i)^{-1} .\]
Because $\|[f(s^dD),h]\|\to 0$ as $s\to 0$, we can write 
\[ \limsup_{s\to 0} \|\,h f(s^dD_m) h - h f(s^dD) h \,\|
=  \limsup_{s\to 0} \;t^d\,\|\, hf(s^dD_m)(D-D_p)hf(s^dD)\,  \| .\]
We decompose the product of operators as follows,
\begin{align*}
  s^d \,&\| hf(s^dD_m)(D-D_m)hf(s^dD)  \| \\
& = s^d \,\|hf(s^dD_m)\|_{L^2\to L^2} \,\|(D-D_m)h\|_{W^d\to L^2} \,\|f(s^dD)\|_{L^2\to W^d}  \\ 
& \le s^d\, \|(D-D_m)h\|_{W^d\to L^2}\,\|f(s^dD)\|_{L^2\to W^d}.
\end{align*}
By Lemma \ref{hypind50:lemm:AftdP} below, the norm $\|f(s^dD)\|_{L^2\to W^d}$ is of order $\OO(s^{-d})$, which means that 
\[ s^d\, \|f(s^dD)\|_{L^2\to W^d} < C.\]
For the remaining factor observe that, because $D_m$ is the principal part of $D$ at $m$,
the norm $\|(D-D_m)h\|_{W^d\to L^2}$ can be made arbitrarily small by choosing a sufficiently small neighborhood of $m$ as the support of $h$.
This proves the last property in Proposition \ref{hypind40:prop:asymp3}, and completes the proof.

\QED

We have to prove a few lemmas that played a role in the previous proof. 
In all the lemmas below, $M$ is a compact contact manifold,
and $D$ is a formally selfadjoint, subelliptic operator of order $d$ with Rockland model operators.
The first of these lemmas is an easy result in the {\em pseudodifferential} Heisenberg calculus.
We give an independent proof.

\begin{lemma}\label{hypind99:lemma:NegativeOrder}
If $k$ is a positive integer then $(D^2+1)^{-k/2d}$ is an operator of order $-k$, in the sense that it is bounded as an operator $L^2\to W^k$. 
\end{lemma}
{\bf Proof.}
We write $\Delta = (D^2+1)^{1/2d}$.
First we show that the operators $\Delta^r$ $(r\in \RR)$ map $C^\infty(M)$ bijectively to $C^\infty(M)$.
Observe that $C^\infty(M)$ is the intersection of the Sobolev spaces $W^k$, 
which is equal to the intersection of the domains of the subelliptic operators $\Delta^{4dk}$ $(k=1,2,3,\ldots)$.
Let $E$ is the projection valued spectral measure of the invertible, selfadjoint operator $\Delta$.
Then $u\in C^\infty(M)$ if and only if
\[ \int |\lambda|^N dE_{u,u} < \infty \]
for all integers $N>0$. Here $E_{u,u}$ denotes the positive Borel measure $E_{u,u}(\omega) = \ang{E(\omega)u,u}$.
If $v=\Delta^r u$ then 
\[ \ang{E(\omega)v,v} = \ang{E(\omega)\Delta^{2r}u,u} = \int_\omega \lambda^{2r} dE_{u,u} ,\]  
which shows that $dE_{v,v} = \lambda^{2r}dE_{u,u}$.
It follows that $u \in C^\infty$ if and only if $v\in C^\infty$.

Now, to prove the lemma, let $A$ be a differential operator of Heisenberg order $k$.
We first show that there is a constant $C>0$, such that, 
\[ \|Au\|  \le C\|(D^2+1)^{k/2d} u \| ,\]
for all $u\in C^\infty(M)$.
Both the selfadjoint operators $(A^*A)^{d}$ and $D^{2k}$ are of order $2dk$,
and $D^{2k}$ is a Rockland operator for which the a priori estimates hold,
\[ \|(A^*A)^{d} u\| \le C(\|D^{2k} u\|^2 + \|u\|^2),\]
for all $u\in W^{2dk}$.
Rewriting these estimates, we derive,
\[ 0 < \ang{ (A^*A)^{2d} u,u} \le C\ang{ (D^{4k}+1)u,u} \le C \ang{ (D^2+1)^{2k}u,u} ,\]
which holds for all $u$ in the domain of the Rockland operator $(D^2+1)^{2k}$ which is $W^{2dk}$. The domain of $(A^*A)^{d}$ certainly contains $W^{2dk}$.

From this we derive,
\[  0 < \ang{A^*A u,u} \le C \ang{(D^2+1)^{k/d}u,u},\]  
for all $u\in C^\infty(M)$, or,
\[ \|Au\| \le C\|(D^2+1)^{k/2d}u .\]
By the first result we may substitute $u= (D^2+1)^{-k/2d}v$ with $v\in C^\infty$, and we get 
\[ \|A(D^2+1)^{-k/2d}v \| \le C\|v\| .\]
for all $v\in C^\infty$, and hence for all $v\in L^2$. 

\hfill{} $\Box$

\noindent{\bf Remark.} We made use of the fact that if $S,T$ are two essentially selfadjoint differential operators with
$0\le \langle Su,u \rangle \le \langle Tu,u\rangle$,
for all $u\in C_c^\infty(M)$,
then for any $0< r < 1$ and $u\in C_c^\infty(M)$
\[ \ang{S^ru,u} \le \ang{T^r u,u}.\] 
This is a version for (unbounded) differential operators of Proposition 1.3.8 in [Pe].
For a careful proof, see [Er].

\begin{lemma}\label{hypind99:lemma:bounds}
Let $r$ be a fixed real number, $0< r \le 1$. 
The supremum of the family of bounded functions $f_t\in C_b(\RR)$ with $t>0$,
\[ f_t(x) = \frac{(x-i)^r}{tx - i},\]
behaves asymptotically as $\|f_t\|_\infty \sim C t^{-r}$ as $t\da 0$. Here $C>0$ is a constant that depends on $r$.
\end{lemma}

\begin{lemma}\label{hypind50:lemm:AftdP}
\[ \|\, (t^dD-i)^{-1} \,\|_{W^k} \le Ct^{-k} .\]
\end{lemma}
{\bf Proof.}
Let $A$ be a differential operator of Heisenberg order $k\le d$. 
We write,
\[ A(t^dP-i)^{-1} = A(P-i)^{-k/d} \;\cdot\;  (P-i)^{k/d} (t^dP - i)^{-1}  .\]
According to Lemma \ref{hypind99:lemma:bounds}, we have the asymptotic behaviour,
\[ \| (P-i)^{k/d} (t^dP - i) \| \le Ct^{-k} .\]
Lemma \ref{hypind99:lemma:NegativeOrder} shows that $A(P-i)^{-k/d}$ is bounded.

\hfill $\Box$

\begin{lemma} \label{hypind99:lemma:commutator}
\[ \lim_{t\to 0} \|\,[(t^dD-i)^{-1},\varphi]\,\| = 0 .\]
\end{lemma} 
{\bf Proof.}
We have,
\[ [(t^dD-i)^{-1},\varphi] = -t^d (t^dD-i)^{-1}[D,\varphi](t^dD-i)^{-1} .\]
The commutator $[D,\varphi]$ is of Heisenberg order $(d-1)$,
so by Lemma \ref{hypind50:lemm:AftdP}, 
\[  \|\,[D,\varphi](t^dD-i)^{-1}\,\| \le Ct^{d-1} .\] 
 
\hfill{} $\Box$

\addcontentsline{toc}{chapter}{References}

\chapter*{References}

\def\item{\vskip2.75pt
plus1.375pt minus.6875pt\noindent\hangindent1em} \hbadness2500
\tolerance 2500
\markboth{References}{References}

\item{[Ar]} V.\ I.\ Arnold,
{\sl Mathematical methods in classical mechanics},
GTM 60, Springer Verlag, Berlin and New York, 1978.


\item{[BG]} R.\ Beals and P.\ Greiner,
{\sl Calculus on Heisenberg Manifolds},
Annals of Mathematics Studies (119), Princeton, 1988.


\item{[Bo]} L.\ Boutet de Monvel, 
{\sl On the index of Toeplitz operators of several complex variables},
Invent. Math. 50 (1979), 249--272.


\item{[CG]} L.\ Corwin and F.P.\ Greenleaf,
{\sl Representations of nilpotent Lie groups and their applications;
Part 1: Basic theory and examples},
Cambridge studies in advanced mathematics 18, 
Cambridge University Press, 1990.

\item{[CGGP]}, M.\ Christ, D.\ Geller, P.\ Glowacki, L.\ Polin,
{\sl Pseudodifferential operators on Groups with dilations},
Duke Math.\ J\. 68 (1992), 31--65.  

\item{[Co1]} A.\ Connes,
{\sl An analogue of the Thom isomorphism for crossed products of a $C^*$-algebra by an action of {\bf R}},
Adv.\ Math.\, 39 (1981), 31--55.

\item{[Co]} A.\ Connes,
{\sl Noncommutative Geometry},
Academic Press, 1994.



\item{[EM]} C.\ Epstein and R.\ Melrose,
{\sl Contact degree and the index of Fourier integral operators},
Math. Res. Lett. 5 (1998), no.3, 363--381.

\item{[EM2]} C.\ Epstein and R.\ Melrose,
{\sl The Heisenberg algebra, index theory and homology},
preprint, 2003.

\item{[Ep]} C.\ Epstein,
{\sl Lectures on Indices and Relative Indices on Contact and CR-manifolds},
Woods Hole Mathematics: Perspectives in Mathematics and Physics,
World Scientific, 2004.

\item{[Er]} E. van Erp,
{\sl The Atiyah-Singer Formula for Subelliptic Operators on a Contact Manifold},
Ph.D. thesis, The Pennsylvania State University, 2005. 

\item{[Er2]} E. van Erp,
{\sl The Atiyah-Singer Formula for Subelliptic Operators on a Contact Manifold, Part II},
to be published.


\item{[FS1]} G.\ B.\ Folland and E.\ M.\ Stein,
{\sl Estimates for the $\bar{\partial}_b$ complex and analysis on the Heisenberg group},
Comm.\ Pure and Appl.\ Math.\, vol XXVII (1974), 429--522.

\item{[FS]} G.\ B.\ Folland and E.\ M.\ Stein,
{\sl Hardy Spaces on Homogeneous Groups},
Princeton University Press, 1982.


\item{[Hi]} N.\ Higson, 
{\sl On the $K$-theory proof of the index theorem},
Index Theory and Operator Algebras (Boulder, CO, 1991), Contemp.\ Math.\, 148, 67--86.


\item{[HN]} B.\ Helffer and J.\ Nourrigat,
{\sl Charact\'erisation des op\'erateurs hypoelliptique homog\`enes invariants a gauche sur un groupe de Lie nilpotent gradu\'e},
Comm.\ Partial.\ Diff.\ Eq.\, 4 (1979), 899--958.






\item{[Ma]} K.\ Mackenzie,
{\sl Lie Groupoids and Lie Algebroids in Differential Geometry},
LNM Lecture Note Series (124), Cambridge University Press, 1987.

\item{[Me]} R.\ Melrose,
{\sl Homology and the Heisenberg algebra}, 
joint work with C.\ Epstein and G.\ Mendoza,
S\'eminaire sur les Equations aux D\'eriv\'ees Partielles, 1996--1997,
Exp. No. XII, \'Ecole Polytech., Palaiseau, 1997.



\item{[Ni2]} V.\ Nistor,
{\sl An index theorem for gauge-invariant families: the case of solvable groups},
Acta Math.\ Hungar.\ 99 (1--2) (2003), 155--183.

\item{[NS]} E.\ Nelson and W.\ F.\ Steinspring,
{\sl Representation of elliptic operators in an enveloping algebra},
Amer.\ J.\ Math.\ 81 (1959), 547--560. 



\item{[Pe]} G.\ Pedersen,
{\sl $C^*$-algebras and their Automorphism Groups},
Academic Press, 1979.










\item{[Ta1]} M.\ E.\ Taylor,
{\sl Noncommutative microlocal analysis, part I},
Mem.\ Amer.\ Math.\ Soc.\, vol. 313, AMS, 1984.

\item{[Ta]} M.\ E.\ Taylor,
{\sl Noncommutative Harmonic Analysis},
AMS, 1986.


\end{document}